\documentstyle[amscd,12pt,psamsfonts]{amsart}
\input xypic

\newtheorem{defn}{Definition}
\newtheorem{thm}[defn]{Theorem}

\newtheorem{lem}[defn]{Lemma}

\theoremstyle{remark}
\newtheorem{rem}{Remark}
\theoremstyle{remark}

\numberwithin{equation}{section}
\numberwithin{defn}{section}

\begin{document}


\newcommand\spk{{\operatorname{Spec}}(k)}
\renewcommand\sp{\operatorname{Spec}}       
\renewcommand\sf{\operatorname{Spf}}       
\newcommand\proj{\operatorname{Proj}}
\newcommand\aut{\operatorname{Aut}}
\newcommand\lie{\operatorname{Lie}}
\newcommand\grv{{\operatorname{Gr}}(V)}
\newcommand\gr{\operatorname{Gr}}
\newcommand\glv{{\operatorname{Gl}}(V)}
\newcommand\glve{{\widetilde{\operatorname{Gl}}(V)}}
\newcommand\gl{\operatorname{Gl}}
\renewcommand\hom{\operatorname{Hom}}
\newcommand\End{\operatorname{End}}
\newcommand\pic{\operatorname{Pic}}
\renewcommand\det{\operatorname{Det}}    
\newcommand\tr{\operatorname{Tr}}      
\newcommand\detd{\operatorname{Det}^\ast}
\newcommand\im{\operatorname{Im}}
\newcommand\res{\operatorname{Res}}
\renewcommand\deg{\operatorname{deg}}
\newcommand\id{\operatorname{I}}
\newcommand\rank{\operatorname{rank}}
\newcommand\limi{\varinjlim}
\newcommand\limil[1]{\underset{#1}\varinjlim\,}
\newcommand\limp{\varprojlim}
\newcommand\limpl[1]{\underset{#1}\varprojlim\,}
\newcommand\dk{{\bf\delta}}       

\renewcommand\o{{\mathcal O}}      
\renewcommand\L{{\mathcal L}}       
\newcommand\B{{\mathcal B}}       
\renewcommand\c{{\widehat C}}
\renewcommand\P{{\mathcal P}}    
\newcommand\I{{\mathcal I}}    
\newcommand\D{{\mathcal D}}    
\newcommand\C{{\mathbb C}}    
\newcommand\Z{{\mathbb Z}}    
\newcommand\M{{\mathcal M}}
\newcommand\T{{\mathcal T}}
\newcommand\U{{\mathcal U}}
\newcommand\W{{\mathcal W}}

\newcommand\w{\widehat}
\renewcommand\tilde{\widetilde}  
\newcommand\f{\overset{\to}f}

\newcommand\iso{@>{\sim}>>}
\renewcommand\lim{\limpl{A\in\B}}
\newcommand\beq{       
      \setcounter{equation}{\value{defn}}\addtocounter{defn}1
      \begin{equation}}

\renewcommand{\thesubsection}{\thesection.\Alph{subsection}}

\title[Generalized KP Hierarchy for Several
Variables]{Generalized KP Hierarchy\\
  for Several Variables}
\author[F. J. Plaza Mart\'{\i}n]{F. J. Plaza Mart\'{\i}n}

\address{ Departamento de
Matem\'aticas, Universidad de Salamanca,  Plaza de la Merced 1-4\\
Salamanca 37008. Spain.}

\thanks{
This work is partially supported by the CICYT research contract n.
PB96-1305 and Castilla y Le\'on regional government contract SA27/98.}

\email{fplaza@@gugu.usal.es}

\begin{abstract}
Following the techniques of M. Sato (see \cite{Sa}), a
generalization of the KP hierarchy for more than one variable is
proposed. An approach to the classification of solutions and a method to
construct algebraic solutions is also offered.
\end{abstract}


\maketitle



\section{Introduction}

After the
remarkable papers of Krichever (\cite{Kr}), Mulase  (\cite{Mulase}),
Segal-Wilson  (\cite{SW}) and Sato (see \cite{Sa} and references therein)
the KP hierarchy was extensively studied and a lot of important results
were eventually given in many related problems. Let us cite some of
their topics: infinitesimal transformations of soliton equations
(\cite{DJKM}); characterization of Jacobian varieties (\cite{Shiota});
relation of commutative rings of differential operators and abelian
varieties  (\cite{Naka}); and more recently, study of the moduli space of
pointed curves  (\cite{MP}); algebraic solutions of the multicomponent KP
hierarchy (\cite{P1}); generalization of the Krichever correspondence for
varieties of dimension greater than
$1$  (\cite{Osipov}); generalization of the KP formalism for
pseudodifferential operators in several variables (\cite{Parshin}); etc. .

This paper aims at generalizing the theory of
the KP hierarchy for several variables following Sato's techniques (see,
for instance,  \cite{Sa}). We think that our approach will be
useful to study the above cited topics in greater generality.

The two main parts in Sato's approach to the KP hierarchy are the use of a
universal grassmann manifold as classifying space for the solutions, on
the one hand, and the construction of solutions from some
algebro-geometric data through the Krichever map, on the other.
Usually, these two points have been studied separately. We offer a
generalization which keeps both aspects together as it was proposed in
\S3.4 of \cite{Sa} for dimension greater than $1$; in this sense, it is
closer to the methods employed in \cite{DJKM,Parshin,Osipov,P1}.

Let us summarize the contents of the paper. The second section begins
with some facts on pseudodifferential operators and the
definition of KP hierarchy for $N$ variables (see, for
instance, \cite{Parshin}). Introducing the notion of wave function for
this hierarchy, one proves that this hierarchy can be written as the
compatibility condition of a system of differential equations (Theorem
\ref{thm:1}). 

In the third section, we study the set of wave functions in terms of an 
infinite dimensional Grassmannian, which is a generalization of those
given in \cite{AMP,P1,Sa,SW}. Theorems \ref{thm:kpgr} and \ref{thm:grkp}
generalize Sato's classification of solutions of the KP hierarchy: there
is a 1-1 correspondence between the set of wave functions and the points
of the infinite Grassmannian.

The last section offers a method for constructing solutions for the new
hierarchy starting with some algebro-geometric data (Theorem
\ref{thm:krich}). Similarly to the standard KP, the solutions obtained
in this way are finite gap. This construction generalizes the Krichever
map. It would be interesting to include Osipov's generalization
(\cite{Osipov}) in our framework.

I would like to express my gratitude to Prof. G. B. Segal for inviting me
to the DPMMS at University of Cambridge (UK) where most of this work has
been done.

\section{The new hierarchy}
\subsection{Pseudodifferential Operators}

Let us begin recalling some standard definitions and
properties of pseudo-differential operators (pdo) for $N$
variables, $N\geq 1$ (see \cite{Parshin}).

Multiindexes will be denoted with greek letters and the following
notation will be used:
\begin{itemize}
\item the entries of a multiindex $\alpha$ will be denoted with
subindexes, i.e. $(\alpha_1,\ldots,\alpha_N)$;
\item $0$ denotes the multiindex $(0,\ldots,0)$;
\item  for $\alpha$ and $\beta$ we say that $\alpha\subseteq \beta$ if
$\alpha_i\leq\beta_i$ for all $i$;
\item  $\alpha\subset \beta$ means that $\alpha\subseteq \beta$ and that
$\alpha \neq \beta$;
\item $\leq$ will denote the reverse lexicographic order in $\Z^N$;
that is, $\alpha\leq \beta$ if $\sum_{i=1}^N \alpha_i k^i\leq
\sum_{i=1}^N \beta_i k^i$ for all $k>>0$.
\end{itemize}

For unknowns  $x_1,\dots,x_N$, let $\C_x:=\C((x_1))\dots((x_N))$ be the
$N$-dimensional local field of iterated Laurent series; that is,
$\C((x_1))\dots((x_i))$ is defined to be the quotient field of
$\C((x_1))\dots((x_{i-1}))[[x_i]]$. Let
$\C[[x]]$ denote $\C[[x_1,\dots,x_N]]$. 

Now, one  considers the $\C[[x]]$-module of pdo's: 
$$\P\,:=\,\left\{\sum a_\alpha\partial^\alpha\,\vert\, a_\alpha\in A, 
n\in\Z\right\}
\subseteq \C[[x]]((\partial_1^{-1}))\dots((\partial_N^{-1}))$$
 where $\partial_i:=\frac{\partial}{\partial x_i}$ ($1\leq i\leq N$) and
$\partial^\alpha:=\partial_1^{\alpha_1}\dots\partial_N^{\alpha_N}$.

Since $[\partial_i,\partial_j]=0$, the following generalization of the
Leibnitz rule: 
$$\big( \sum_{\alpha} a_{\alpha}\partial^{\alpha}\big)\big( \sum_{\beta}
b_{\beta}\partial^{\beta}\big)
\,:=\, 
\sum_{\alpha,\beta}\sum_{0\subseteq\gamma}\prod_{i=1}^N\binom{\alpha_i}{\gamma_i}a_\alpha(\partial^\gamma
b_\beta)\partial^{\alpha+\beta-\gamma}$$ 
endows $\P$ with a ${\mathbb C}$-algebra structure. Mapping a pdo
$P=\sum_\alpha a_\alpha\partial^\alpha$ to $ P_+:= \sum_{\alpha\geq
0}a_\alpha\partial^\alpha$ one obtains an endomorphism $\P\to\P$ whose
image (resp. kernel) will be denoted by $\P_+$ (resp. $\P_-$).

\subsection{KP for $N$ variables}

Let $\C[[t]]$ denote $\C[[\{t_\alpha\}_{0\subset\alpha}]]$ and
$\partial_{\alpha}$ be $\frac{\partial}{\partial t_\alpha}$. For pdo
$L_1,\dots, , L_N$ let 
$L^\alpha$ denote $L_1^{\alpha_1}\cdots L_N^{\alpha_N}$.

\begin{defn}
The KP hierarchy for $N$ variables, KP(N), is the
following Lax system:
\beq
\cases
\partial_{\alpha}L_i\,=\,\lbrack(L^\alpha)_+,L_i \rbrack 
\\
\lbrack L_i,L_j\rbrack \,=\, 0  
\endcases
\qquad 1\leq i,j \leq N\, , \, 0\subset\alpha
\label{eq:KPN-Lax}
\end{equation}
where:
\beq
L_i\,=\, \partial_i+\sum_{\alpha\subset 0}
u_{i\alpha}(t)\partial^\alpha\in \P\otimes\C[[t]]
\qquad i=1,\dots,N
\label{eq:type}\end{equation}
\end{defn}

A formal oscillating function, $w(t,x)$, over this ring is a formal
expression of the following type:
$$\big(1+\sum_{0\subset\alpha}
a_\alpha(t) x^\alpha\big)\cdot e^{\xi(t, x)}$$
where $e^{\xi(t, x)}:=\exp(\sum_{0\subset\beta} t_\beta x^{-\beta})$,
$t=\{t_\beta\,\vert\,0\subset\beta\}$ are ``time'' variables and 
$a_\alpha(t)\in\C[[\{t_\alpha\}_{0\subset\alpha}]]$.

The KP hierarchy as defined in \ref{eq:KPN-Lax} consists of a set of
deformation equations for the operators $L_i$. Similarly to \S1.2 of
\cite{Shiota}, it admits other equivalent formulations.

\begin{thm}\label{thm:1}
The following conditions are equivalent:
\begin{enumerate}
\item $\{L_1,\dots,L_N\}$ are pdo of the type
\ref{eq:type} and satisfy \ref{eq:KPN-Lax};
\item there exists a pdo $S\in 1+\P_-\otimes\C[[t]]$ such that:
\beq
\cases L_i S\,=\, S\partial_i \qquad  &  1\leq i\leq N
\\
\partial_\alpha S\,=\, -(S\partial^\alpha S^{-1})_-S\qquad &0\subset
\alpha
\endcases
\label{eq:KPNpdo}\end{equation}
\item $\{L_1,\dots,L_N\}$ admits a wave function; that is, a
formal oscillating function $w(t,x)$ such that:
\beq
\cases
L_i w\,=\, x_i\cdot w\qquad & i=1,\dots,N \\
\partial_\alpha w\,=\, (L^\alpha)_+ w \qquad&  0\subset\alpha
\endcases
\label{eq:KPNwave}\end{equation}
\end{enumerate}
\end{thm}

\begin{pf}
The equivalence of 2 and 3 is formal since $w(t,x)$ and
$S(t,\partial)$ are related through $w(t,x)=S(t,x)e^{\xi(t,x)}$.

From \cite{Parshin} one knows that 2 implies 1. Conversely, if
$\{L_1,\dots,L_N\}$ are of the form \ref{eq:type} and commute pairwise,
then there exists $S\in 1+\P_-$ such that $L_i=S\partial_i S^{-1}$
(\cite{Parshin}), and such $S$ is unique up to right multiplication
by a pdo in $1+\P_-$ with constant coefficients. Further, arguments
similar to those of \S1.2 of \cite{Shiota} show that equation
\ref{eq:KPNpdo} determines recursively the coefficients of
$S$.
\end{pf}

\begin{defn}
A finite gap wave function for the KP(N) hierarchy is a wave function
$w(t,x)$ such that $\partial_\alpha w=0$ for a finite number of
$\alpha$'s.
\end{defn}

\section{Classification of Solutions}

From now on $V$ will denote the complex vector space $\C_x$ endowed with
the filtration
$V_n:=x_N^n\C((x_1))\dots((x_{N-1}))[[x_N]]$ with $n\in\Z$.

For a subspace $A$ of $V$, one observes that:
$$\begin{aligned}
d_A: \Z &\to \Z\cup\{\infty\}  \\
n &\mapsto \dim(V_n/A\cap V_n)
\end{aligned}$$
is a decreasing function. 

\begin{defn}
We consider the linear topology in $V$ given by the
basis, $\B$, of neighbourhoods of $(0)$: the set of proper subspaces
$A\subset V$ such that $d_A(n)$ is finite for all $n$ and converges to
zero.
\end{defn}

Given $A\in\B$ and a proper subspace of $V$, $B$, the
following properties hold:
\begin{itemize}
\item if $A\subseteq B$, then $B\in\B$;
\item if $B\subseteq A$ is of finite codimension, then $B\in\B$;
\item if $B\in \B$, then $A\cap B\in\B$.
\end{itemize}

\begin{defn}
The Grassmannian of the pair $(V,\B)$ is the $\C$-scheme representing
the functor:
$$S\rightsquigarrow \left\{\begin{gathered}
\text{submodules }U\subset \hat
V_S\text{ such that }U\oplus \hat A_S\to\hat V_S  \\
\text{ is an isomorphism for a subspace }A\in\B
\end{gathered}\right\}$$ 
(where $\hat A_S$ denotes the completion of $A\otimes_{\C}\o_S$).
\end{defn}

Recall that in \cite{AMP,MP,P2} the Grassmannian
of $(V=\C((x)), V_n=x^n\C[[x]])$ is defined using the notion of
subspace commensurable with $V_0$ and observe that such notion  is
more restrictive. However, the existence of the scheme
$\grv$ is, in both cases, a direct consequence of the following:

\begin{lem}
Let $A$ be an element of $\B$. Then, $\hat A$ is an inverse limit of
finite dimensional $\C$-vector spaces. 
\end{lem}

\begin{pf}
Note that $\B$ gives a basis of the induced topology on $A$;
namely $\{A\cap B\vert B\in\B\}$. 

For a given subspace $B\in\B$, the sets:
$$\begin{aligned}
\B_{A}\,&:=\,\{C\subseteq A \text{ s.t. }\dim
A/C<\infty\}
\\
\B_{A,B}\,& :=\, \{C\in \B_{A}
\text{ s.t. }A\cap B\subseteq C\subseteq A\}
\end{aligned}
$$
satisfy $\B_{A,B}\subseteq
\B_A\subseteq \B$. The commutativity of the
following diagramm:
$$
\xymatrix{
\hat A=\limpl{B\in\B} A/(A\cap B) \ar@{->>}[r] \ar@{->>}[d] & 
\limpl{C\in\B_{A}} A/(A\cap C) \ar@{->>}[d] \\
A/(A\cap B) \ar@{^{(}->}[r] & \limpl{C\in\B_{A,B}} A/(A\cap C)}
$$
implies that the horizontal arrow on the top is an isomorphism, and the
conclusion follows.
\end{pf}

Motivated by the properties of Baker-Ahkiezer functions for the standard KP
hierarchy (see \cite{Kr,MP,SW}), we give the following: 

\begin{defn}
Given a point $U\in\grv$ and a basis $\{u_i\}_{i\in I}$ of it, 
define the Baker-Ahkiezer function,
$\omega_U$, as the formal sum:
$$\omega_U(t,x)\,:=\, \sum_{i\in I} t_{v(u_i)} u_i$$
where $v:V\to \Z^N$ maps $\sum_\alpha
a_\alpha x^\alpha$ to the multiindex $\alpha$ such that $a_\beta=0$ for
all $\beta<\alpha$.
\end{defn}


Define $U_0$ to be the subspace $\C[x_1^{-1},\dots,x_N^{-1}]$ and
consider $A_0\subset V$ satisfying:  $U_0\oplus A_0\simeq V$;
$x^\alpha\in A_0$ for all multiindex $\alpha\not\subset 0$; and
$A_0\in \B$. Let $F_0$ be the open subscheme $\hom(U_0,\hat
A_0)\subset\grv$.

\begin{thm}\label{thm:kpgr}
There is a natural map from the set of wave functions for the
KP(N) hierarchy to $F_0\subset \grv$.
\end{thm}

\begin{pf}
Let $w(t,x)=S(t,\partial_t)e^{\xi(t,x)}$ be a formal oscillating
function with $S(t,x)=1+\sum_{\beta<0} a_\beta x^\beta$. 

If $w$ is a wave function for the KP(N) hierarchy and $0\subseteq\alpha$ is
a multiindex, one has:
$$\partial_\alpha w\,=\, x^{-\alpha} w- (L^\alpha)_-w\,=\,
\big(x^{-\alpha} S-(S\partial^\alpha S^{-1})_-S\big)e^{\xi(t,x)}$$
and therefore:
$$(\partial_\alpha w)\vert_{t=0}\,=\, x^{-\alpha}+ (\text{higher order
terms})$$

This means that $(\partial_\alpha w)\vert_{t=0}$ generates a point of
$F_0\subset\grv$ as $\alpha$ varies.
\end{pf}

\begin{rem}
 This map is equivariant w.r.t. the
action by right multiplication of
$1+\C((\partial_1^{-1}))\dots((\partial_N^{-1}))_-$ on the set of wave
functions and that of $1+\C((x_1^{-1}))\dots((x_N^{-1}))_-$ on
$\grv$.
\end{rem}

\begin{thm}\label{thm:grkp}
If $U$ is a point of $F_0$, then its Baker-Ahkiezer function is a wave
function for the KP(N).
\end{thm}

\begin{pf}
Let $K$ be the field $\C((x_1))\dots((x_{N-1}))$. If $U$ is a subspace 
in $F_0$, then the $K$-vector space $\tilde U\subset K((x_N))$ given by
the image of $K\otimes_{\C} U\to K((x_N))$ is a point of
$\gr(K((x_N)),K[[x_N]])$ (see \cite{AMP,MP} for the definition of the
Grassmannian). Since  $K((x_N))$ carries a non-degenerate bilinear pairing,
$<\, ,\,>$, induced by the residue at $x_N=0$, we consider $(\tilde
U)^\perp\in \gr(K((x_N)),K[[x_N]])$. Let
$\psi_U(s,x_N)$ be the wave function of $(\tilde U)^\perp$.

From the definition of $\omega_U(t,x)$ and from the properties of
$\psi_U(s,x_N)$, one obtains that:
$$<\omega_U(t,x),\psi_U(s,x_N)>\,=\,0\qquad \forall t,s$$

This identity implies that (see \cite{DJKM}):
$$\big((\partial_\alpha+(L^\alpha)_-)
S(t,\partial_t) R^*(t,\partial_N)\big)_-\,=\,0 \qquad\forall \alpha$$
where $\omega_U(t,x)=S(t,\partial_t) e^{\xi(t,x)}$, 
$\psi_U(s,x_N)=R(s,\partial_s)e^{\sum_{i>0}s_ix_N^{-i}}$, and $R^*$ is
the adjoint operator of $R$. Then, it follows that:
$$(\partial_\alpha+(L^\alpha)_-) S(t,\partial_t)\,=\,0 \qquad\forall
\alpha$$
\end{pf}

\begin{rem}
Once one has related the KP(N) hierarchy with the infinite Grassmannian
$\grv$, one observes that the groups acting on it will give
symmetries of the hierarchy (see \cite{DJKM} for the case of the KP
hierarchy).

Further, there is a natural 2-cocycle given by 
$det(\delta_B^{-1}\circ\delta_A)\in H^0(F_A\cap F_B,\o^*_{F_A\cap F_B})$,
where $\L$ is the universal submodule
and $\delta_A$ is the morphism $\L\oplus A\to V$. Then, for a group acting
on $\grv$ the line bundle  associated to this cocycle defines central
extensions of the group and of its Lie algebra.

The case of the Lie algebra of the group $\aut_{\C-\text{alg}}\hat V$
has beautiful properties when $N=1$ due to its relation with the
Virasoro algebra (see
\cite{MP2}). However, the
$N>1$ situation differs substantially since the Lie algebra has no
non-trivial central extensions. This fact follows from \cite{fisica}, since
that Lie algebra of
$\aut_{\C-\text{alg}}\hat V$ has a dense subalgebra with
generators $\{L^a_\alpha\,\vert
\alpha\in\Z^N\, ,\, 1\leq a\leq N\}$
and the following Lie bracket $[L^a_\alpha,L^b_\beta]\,=\, \beta_a
L^b_{\alpha+\beta}-\alpha_b L^a_{\alpha+\beta}$.
\end{rem}

 \section{Finite Gap Solutions}

From \cite{Kr,SW} it is known that for the standard KP hierarchy the
Krichever morphism provides a way to construct solutions starting with the
the geometric data $(C,p,\alpha)$ where $C$ is an algebraic
integral curve over $\C$, $p\in C$ is a smooth point, and $\alpha$ is
an isomorphism $\hat\o_{C,p}\simeq \C[[x_0]]$.

Now, we will see how finite gap solutions for the KP(N) hierarchy 
might be constructed from some algebro-geometric data. 

\begin{defn}
An algebro-geometric datum for the KP(N) consists of
$(X,p,\alpha,\{Y_1,\dots,Y_N\})$ where $X$ is a
$N$-dimensional integral regular projective scheme, $p$ is a
point, $\alpha$ is an isomorphism $\hat\o_{X,p}\simeq
\C[[x_1,\dots,x_N]]$, and
$\{Y_1,\dots,Y_N\}$ is a ordered set of Weil divisors such that:
$p\in\cap_i Y_i$, and
 $x_i=0$ is the local equation of $Y_i$ in a
neighbourhood of $p$.
\end{defn}

Observe that given an algebro-geometric datum
$(X,p,\alpha,\{Y_1,\dots,Y_N\})$, the morphism $\o_{X,p}
\hookrightarrow\hat\o_{X,p}\iso \C[[x_1,\dots,x_N]]$ induces a map from
the function field of $X$, $\Sigma$, to $V$. Moreover, given a function
$f\in\Sigma$, one can define $v(f)\in\Z^N$
as the smallest (w.r.t. the order
$\leq$) exponent occurring in the image of $f$ by the natural map
$\Sigma\to V$.

Now, the map $\Sigma\hookrightarrow V$ and Serre's vanishing
theorem allow an easy generalization of Krichever's construction
(\cite{Kr} in the following form:
 
\begin{thm}\label{thm:krich}
Let ${\mathfrak X}=(X,p,\alpha,\{Y_1,\dots,Y_N\})$ be an algebro-geometric
datum for the KP(N). Then, the $\C$-vector space:
$$A_{\mathfrak X}:=\{f\in \underset{0\subseteq
\alpha}\bigcup H^0(X,\o_X(\sum_i\alpha_i
Y_i)\,\vert\, 0\subseteq v(f)\}$$ is a point of $\grv$ and its wave
function is a finite gap solution of the KP(N) hierarchy.
\end{thm}

To finish this section let us point out two remarks. First, it could be
interesting to adapt Osipov's generalization of the Krichever
correspondence (\cite{Osipov}) to this kind of algebro-geometric data.
Second, one wonders about the possible generalizations of the
Burchnall-Chaundy theory for studying rings of commuting differential
operators within this framework.



\end{document}